\documentclass[12pt,a4paper,dvips]{article}
\usepackage[cp1251]{inputenc}
\usepackage[english]{babel}
\usepackage{amsmath,amssymb,amsthm}
\usepackage[dvips,pdftex]{graphicx}
\usepackage{multirow}
\usepackage{rotating}
\usepackage{geometry,comment}
\usepackage[dvips,pdftex]{xcolor}
\usepackage{authblk}
\usepackage{epstopdf}

\newcommand{\SO}{\mathrm{SO}}

\newcommand{\SU}{\mathrm{SU}}

\newcommand{\su}{\mathfrak{su}}

\newcommand{\id}{\mathrm{id}}

\newcommand{\p}{\bar{p}_3}
\newcommand{\st}{\sin{\tau}}
\newcommand{\ct}{\cos{\tau}}

\newcommand{\se}{\sin{(\tau \eta \p)}}
\newcommand{\ce}{\cos{(\tau \eta \p)}}

\def\tcut{t_{\operatorname{cut}}}
\def\taucut{\tau_{\operatorname{cut}}}
\def\tauconj{\tau_{\operatorname{conj}}}

\DeclareMathOperator{\diam}{diam}

\theoremstyle{definition}

\newtheorem{Remark}{Remark}
\newtheorem*{Acknowledgement}{Acknowledgement}

\theoremstyle{plain}

\newtheorem{Theorem}{Theorem}

\newcommand*{\affaddr}[1]{#1} 

\newcommand*{\email}[1]{\texttt{#1}}

\title{Diameter of $\SU_2$ for a left-invariant axisymmetric Riemannian metric
\footnote{This work is supported by the Russian Science Foundation
under grant 17-11-01387 and performed in Ailamazyan Program Systems
Institute of Russian Academy of Sciences.}
}

\author{A.~V.~Podobryaev\\
\affaddr{Program Systems Institute of RAS}\\
\email{alex@alex.botik.ru}
}

\begin{document}

\maketitle

\begin{abstract}
We consider the Lie group $\SU_2$ endowed with a left-invariant axisymmetric Riemannian metric.
This means that a metric has eigenvalues $I_1 = I_2, I_3 > 0$. 
We give an explicit formula for the diameter of such metric. Other words, we compute the diameter of Berger's sphere.



\end{abstract}

\section{\label{section-introduction}Introduction}

Several aspects of analysis of wave, heat, Schr\"{o}dinger equations on a Riemannian manifold $(M, g)$
require the spectral analysis of the Laplace-Beltrami operator $\triangle_g$ that corresponds to the metric $g$.
The bounds and asymptotic for non-zero elements of the $\triangle_g$-spectrum depend on geometrical properties of the manifold $M$,
such as the dimension, the curvature, the volume and the diameter $\diam_g{M}$ of $M$ (see, for example~\cite{lablee}).

Recall that $\diam_g{M} = \sup{\{d_g(x, y) \ | \ x, y \in M\}}$, where $d_g$ is the Riemannian distance.

Denote by $g(I_1, I_2, I_3)$ the left-invariant Riemannian metric on $\SU_2$ with eigenvalues $I_1 \leqslant I_2 \leqslant I_3$.
N.~Eldredge, M.~Gordina and L.~Saloff-Coste~\cite{eldredge-gordina-saloff-coste} show that
\begin{equation}
\label{diam-inequality}
D_0 \sqrt{I_2} \leqslant \diam_{g(I_1, I_2, I_3)}{\SU_2} \leqslant D_{\infty} \sqrt{I_2},
\end{equation}
where $D_0, D_{\infty}$ are some constants.

In this paper we prove an explicit formula of the diameter $\diam_{g(I_1, I_1, I_3)}{\SU_2}$ in the case $I_1 = I_2$.
Note, that $\SU_2$ endowed with the left-invariant metric $g(I_1, I_2, I_3)$ is known as Berger's sphere~\cite{berger}.

{\Theorem
\label{th-diameter}
The diameter of $\SU_2$ for left-invariant Riemannian metric with eigenvalues $I_1 = I_2, I_3 > 0$ is equal to
$$
\diam_{g(I_1, I_1, I_3)}{\SU_2} =
\begin{aligned}
\left\{
\begin{array}{lll}
2\pi\sqrt{I_1},                   & \text{for} & I_1 \leqslant I_3, \\
2\pi\sqrt{I_3},                   & \text{for} & I_3 < I_1 \leqslant 2I_3, \\
\frac{\pi I_1}{\sqrt{I_1 - I_3}}, & \text{for} & 2I_3 < I_1.
\end{array}
\right.
\end{aligned}
$$
}
\medskip

{\Remark
\label{rem-diameter-continuous}
The diameter is a continuous function of the variables $I_1, I_3$.
Indeed, for $I_1 = I_3$, we have $2\pi \sqrt{I_1} = 2 \pi \sqrt{I_3}$.
For $I_1 = 2I_3$, we obtain $\frac{\pi I_1}{\sqrt{I_1 - I_3}} = \frac{2\pi I_3}{\sqrt{I_3}} = 2\pi\sqrt{I_3}$.
}
\medskip

{\Remark
\label{rem-inequality}
The diameter satisfies inequality~(\ref{diam-inequality}).
Indeed, one can put $D_0 = \pi$ and $D_{\infty} = 2\pi$.
}
\medskip

For proof of Theorem~\ref{th-diameter} see Section~\ref{section-proof}.
The proof is based on our previous results~\cite{podobryaev-sachkov} for the cut locus and cut time for axisymmetric left-invariant Riemannian metrics on $\SU_2$.
(Note, that the diameter $\diam_{g(I_1, I_1, I_3)}{\SO_3}$ is computed in that paper.)
Also we use the equation for the conjugate time achieved by L.~Bates and F.~Fass\`{o}~\cite{bates-fasso}.
Section~\ref{section-cut-time} contains a summary of necessary notation and results from paper~\cite{podobryaev-sachkov}.

\section{\label{section-cut-time}Cut time}

We use the Hamiltonian formalism.
Every geodesic (starting form $\id \in \SU_2$) with arc\-length parametrization is defined by an initial momentum $p \in \su_2^*$ such that $H(p) = \frac{1}{2}$,
where $H(p) = \frac{p_1^2}{I_1} + \frac{p_2^2}{I_2} + \frac{p_3^2}{I_3}$ is the Hamiltonian,
and $p_1, p_2, p_3$ are components of $p$ in the basis dual to the basis where the Killing form and the metric are diagonal.
(see~\cite{podobryaev-sachkov}, Section~4).
Introduce the following notation:
$$
|p| = \sqrt{p_1^2 + p_2^2 + p_3^2}, \qquad \p = \frac{p_3}{|p|}.
$$

The cut time $\tcut{(p)}$ is a time of loss of optimality for the geodesic with the initial momentum $p$.
Due to the axisymmetry of the metric, the cut time is a function $\tcut{(\p)}$ of variable $\p \in [-1, 1]$.

Put $\tau = \frac{t|p|}{2I_1}$. Define $\tau_3(\p)$ as the first positive root of the equation
$$
\ct \se + \p \st \ce = 0, \qquad \eta = \frac{I_1}{I_3} - 1.
$$
The function $\tau_3(\p)$ is defined on the domain $[-1, 1] \setminus \{0\}$.
There exists $\lim_{\p \rightarrow 0}{\tau_3(\p)}$ equal to $\tauconj{(0)}$,
where $\frac{2I_1\tauconj{(\p)}}{|p|}$ is a conjugate time for the geodesic with the initial momentum $p$
(see~\cite{podobryaev-sachkov}, proof of Proposition~9).
Define $\tau_3(0) = \tauconj{(0)}$.

{\Theorem
\label{th-summary}
\textup{(1) (\cite{podobryaev-sachkov}, Propositions~8, 10)}
For the cut time $\tcut{(\p)} = \frac{2I_1\taucut{(\p)}}{|p|}$, we have
$$
\taucut{(\p)} =
\begin{aligned}
\left\{
\begin{array}{lll}
\pi,        & \text{for} & \eta \leqslant 0, \\
\tau_3(\p), & \text{for} & \eta > 0.
\end{array}
\right.
\end{aligned}
$$
\textup{(2) (\cite{podobryaev-sachkov}, proof of Proposition~9)}
The function $\tau_3(\p)$ is smooth and increasing at the interval $[-1, 0]$ and decreasing at the interval $[0, 1]$.\\
\textup{(3) (\cite{bates-fasso}, Lemma~5)}
If $\eta > 0$, then $\tauconj{(\p)}$ is the first positive root of the equation
$$
\tan{\tau} = -\tau\eta\frac{1-\p^2}{1 + \eta\p^2}.
$$
The function $\tauconj{(\p)}$ is continuous and $\frac{\pi}{2} < \tauconj{(\p)} \leqslant \pi$.\\
\textup{(4) (\cite{podobryaev-sachkov}, proof of Proposition~10)}
The inequality $\tau_3(\p) < \tauconj{(\p)}$ is satisfied for $\p \neq 0$.
}
\medskip

\section{\label{section-proof}Proof of Theorem~\ref{th-diameter}}

The diameter is equal to the maximum value of the cut time $\tcut{(\p)}$, which is a function of the variable $\p \in [-1, 1]$.

Consider first the case $\eta \leqslant 0$.
It is easy to see that $|p| = \sqrt{\frac{I_1}{1 + \eta\p^2}}$ (see~\cite{podobryaev-sachkov}, proof of Theorem~4).
It follows from Theorem~\ref{th-summary}~(1), that $\tcut{(\p)} = 2\pi\sqrt{I_1}\sqrt{1 + \eta\p^2}$.
This function has a maximum at the point $\p = 0$.
The maximum value is $2\pi\sqrt{I_1}$.

Consider now the case $\eta > 0$.
We find critical points of the function (Theorem~\ref{th-summary}~(1))
$$
\tcut{(\p)} = \frac{2I_1\tau_3(\p)}{|p|} = 2\sqrt{I_1}\tau_3(\p)\sqrt{1+\eta\p^2}.
$$
The function $\tau_3(\p)$ is even (due to the definition), so the function $\tcut{(\p)}$ is even as well.
Hence, we will consider $\p \in [0, 1]$.

Calculate the derivative
$$
\frac{d \tcut}{d \p}(\p) = 2\sqrt{I_1}
\left(
\frac{d \tau_3(\p)}{d \p}\sqrt{1+\eta\p^2} + \frac{\tau_3(\p) \eta \p}{\sqrt{1+\eta\p^2}}
\right).
$$
To make formulas more compact, we will omit the argument $\p$ of the function $\tau_3(\p)$ below.
After transformations, using the formula for $\frac{d \tau_3}{d \p}$ (see~\cite{podobryaev-sachkov}, proof of Proposition~9)
$$
\frac{d \tau_3}{d \p} =
-\frac{\tau_3\eta\cos{\tau_3}\cos{(\tau_3\eta\p)} + \sin{\tau_3}\cos{(\tau_3\eta\p)} - \tau_3\eta\p\sin{\tau_3}\sin{(\tau_3\eta\p)}}
{-(1+\eta\p^2)\sin{\tau_3}\sin{(\tau_3\eta\p)} + \p(1+\eta)\cos{\tau_3}\cos{(\tau_3\eta\p)}},
$$
we have (up to a positive multiplier $C$)
$$
\frac{d \tcut}{d \p}(\p) =
C \frac{\cos{(\tau_3\eta\p)}
[-(1+\eta\p^2)\sin{\tau_3} - \tau_3\eta(1-\p^2)\cos{\tau_3}]}
{-(1+\eta\p^2)\sin{\tau_3}\sin{(\tau_3\eta\p)} + \p(1+\eta)\cos{\tau_3}\cos{(\tau_3\eta\p)}}.
$$

Consider the case $\cos{(\tau_3\eta\p)} = 0$.
It follows from the definition of $\tau_3$ that in this case $\cos{\tau_3} = 0$.
So, $\tau_3 = \frac{\pi}{2}$, and $\p = \frac{2k+1}{\eta}$, where $k \in \mathbb{Z}$.
Since $\tau_3$ is decreasing at the interval $[0, 1]$ (Theorem~\ref{th-summary}~(2)),
then at most one point of this series is inside of the interval $[0, 1]$.
There is exactly one point if and only if $\eta \geqslant 1$.
This critical point is $\p = \frac{1}{\eta}$.

If $\cos{(\tau_3\eta\p)} \neq 0$ and $\p \neq 0$, then $\cos{\tau_3} \neq 0$ (due to the definition of $\tau_3$).
Divide the numerator and the denominator by $\cos{\tau_3}\cos{(\tau_3\eta\p)}$.
Then, we have (up to a positive multiplier)
$$
\frac{d \tcut}{d \p}(\p) = \frac{-(1+\eta\p^2)\tan{\tau_3}-\tau_3\eta(1-\p^2)}
{-(1+\eta\p^2)\tan{\tau_3}\tan{(\tau_3\eta\p)} + \p(1+\eta)}.
$$
Due to the definition of $\tau_3$, the denominator is equal to
$$
(1+\eta\p^2)\p\tan^2{\tau_3} + \p(1+\eta) > 0.
$$
The first positive root of the numerator is $\tauconj{(\p)}$ (Theorem~\ref{th-summary}~(3)) and  $\tauconj{(\p)} > \tau_3(\bar{p}_3)$ (Theorem~\ref{th-summary}~(4)).
So, the nominator does not vanish.

If $\p = 0$, then $\tauconj{(0)} = \tau_3(0)$. So, $\p = 0$ is a critical point of $\tcut{(\p)}$.
For small enough $\p > 0$, we have $\tau_3(\p) \in (\frac{\pi}{2}, \tauconj{(\p)})$
(since, $\tau_3(0) = \tauconj{(0)} > \frac{\pi}{2}$ and the function $\tau_3$ is continuous and decreasing).
So, the numerator of $\frac{d \tcut}{d \p}(\p)$ is positive.
Therefore, $\p = 0$ is the minimum point for $\tcut{(\p)}$.

Now we prove that $\p = \frac{1}{\eta}$ is a maximum point for $\eta > 1$.
Compute $\frac{d \tcut}{d \p}(1)$. Up to a positive multiplier, this is equal to
$$
\frac{-(1+\eta)\tan{\tau_3(1)}}{(1+\eta)\tan^2{\tau_3} + (1+\eta)}.
$$
The denominator is positive. One can find $\tau_3(1)$ from the equation
$$
\cos{\tau_3(1)}\sin{(\tau_3(1)\eta)} + \sin{\tau_3(1)}\cos{(\tau_3(1)\eta)} = \sin{(\tau_3(1)(1+\eta))} = 0.
$$
We have $0 < \tau_3(1) = \frac{\pi}{1+\eta} < \frac{\pi}{2}$. Thus, $\tan{\tau_3(1)} > 0$.
So, the function $\tcut{(\p)}$ has two critical points $\p = 0$ and $\p = \frac{1}{\eta}$.
The derivative of this function is positive at the interval $(0, \frac{1}{\eta})$ and negative at the interval $(\frac{1}{\eta}, 1]$.

Finally, if $\eta > 1$, then $\p = \frac{1}{\eta}$ is a maximum point for $\tcut{(\p)}$.
The corresponding value is $\tcut{(\frac{1}{\eta})} = \frac{2 I_1 \tau_3(\frac{1}{\eta})}{|p|}$,
where $\tau_3(\frac{1}{\eta}) = \frac{\pi}{2}$, and $|p| = \sqrt{\frac{I_1}{1 + \frac{1}{\eta}}}$. The diameter is $\pi \sqrt{I_1} \sqrt{1 + \frac{1}{\eta}}$.

If $0 < \eta \leqslant 1$, then the maximum value of $\tcut(\bar{p}_3)$ is achieved at $\bar{p}_3 = 1$. Since $\tau_3(1) = \frac{\pi}{1+\eta}$ and $|p| = \sqrt{\frac{I_1}{1+\eta}}$, the diameter equals $\frac{2 \pi \sqrt{I_1}}{\sqrt{1+\eta}}$.

Using the definition of $\eta$, we obtain the statement of Theorem~\ref{th-diameter}.
$\Box$
\medskip

{\Acknowledgement
We thank E.~Lauret for a question that inspired this work.
}

\end{document}